\documentclass[a4paper]{article}
\pdfoutput=1
\usepackage{amsmath,amsthm,amssymb}
\usepackage{graphicx}
\usepackage[squaren,thickqspace]{SIunits}

\newcommand{\R}{\mathbb{R}}
\newcommand{\Z}{\mathbb{Z}}

\renewcommand{\>}{\rangle}
\newcommand{\ftrue}{f^{\rm true}}
\newcommand{\frec}{f^{\rm rec}}
\newcommand{\gnoise}{g^{\rm noise}}
\newcommand{\Gnoise}{G^{\rm noise}}
\newcommand{\gdata}{g^{\rm data}}
\newcommand{\E}{\operatorname{E}}
\newcommand{\Var}{\operatorname{Var}}
\newcommand{\argmin}{\operatorname{arg\,min}}
\newcommand{\smin}{s_{\rm min}}
\newcommand{\erfc}{\operatorname{erfc}}
\newcommand{\erfcinv}{\erfc^{-1}}

\newtheorem{lemma}{Lemma}

\theoremstyle{definition}

\theoremstyle{remark}
\newtheorem{remark}{Remark}

\title{A new principle for choosing regularization parameter in certain inverse problems}

\author{Hans Rullg\aa rd}

\date{April 15, 2008}

\begin{document}

\maketitle

\begin{abstract}
A new parameter choice rule for inverse problems is introduced. This parameter choice
rule was developed for total variation regularization in electron tomography and might
in general be useful for $L^1$ regularization of inverse problems with high levels of
noise in the data.
\end{abstract}

\section{Introduction}\label{sec:intro}

A standard procedure when solving an ill-posed inverse problem $Tf = g$, where
$T: X \to Y$ is an operator which lacks a continuous inverse, is to define a class of 
continuous operators $S_\lambda: Y \to X$ depending on a parameter $\lambda$ and 
approximating an inverse of $T$ as $\lambda \to 0$. An approximate solution of the
inverse problem is then given by $S_\lambda(g)$ where the regularization parameter
$\lambda$ must be selected by some parameter choice rule. Numerous parameter choice
rules have been proposed, which are suitable for different inverse problems. For a 
nice survey, see \cite[Chapter 4]{EnHaNe96}.

The parameter choice rule proposed in this paper is inspired by the application of
$L^1$-type regularization methods, particularly total variation (TV) regularization, 
to the inverse problem in electron tomography (ET).
For a survey of total variation image restoration methods, see \cite{ChEsPaYi06},
and for a comprehensive account of the mathematics in ET, see \cite{FaOk08}.
It turns out that the well-known parameter choice rules are difficult to apply to this 
problem. One reason for this is that the regularized inverses $S_\lambda$ are non-linear
even if the forward operator $T$ is linear. A second reason is that this inverse problem
has unususally high noise level in the data, the norm of the noise component by far
exceeding the norm of the signal. For this reason, in order to recover anything at all,
it is necessary to make use of knowledge about the statistical properties of the
noise, not just its magnitude.

These issues are taken into account by the proposed parameter choice rule. It is specifically
designed for regularization functionals of $L^1$ type, and for inverse problems aiming at 
reconstructing sparse objects. As it is applied here to ET data, the method is a heuristic,
or error free, parameter choice rule in the terminology of \cite{EnHaNe96}. 
By this one understands a parameter choice rule which 
does not require an explicit estimate of the noise level to be made. Instead, the
noise level is estimated directly from the data, relying on certain assumptions about
the nature of the noise. 

It turns out that the method developed for ET can be 
broken down in several components, which can probably be applied independently in various
other inverse problems. Here an attempt is made to present these components independently.
The paper is organized as follows. In Section~\ref{sec:basic}, the most general part of 
the method, requiring the least structure of the inverse problem, is described.
In Section~\ref{sec:TV} the framework is applied to total variation regularization.
Finally, in Section~\ref{sec:ET}, the application to ET is developed, and numerical 
examples are presented. A general discussion is given in Section~\ref{sec:discussion}.

\section{Basic principle}\label{sec:basic}

Consider the following inverse problem. Let $T: X \to Y$ be a linear operator 
between two linear spaces, with $Y$ a Hilbert space. Let $\ftrue \in X$ be an 
unknown element, which it is our goal to estimate. What is known is an element
$\gdata = T\ftrue + \gnoise \in Y$ where $\gnoise$ is a sample of a random vector
$\Gnoise$. Let us assume that $\E[\Gnoise] = 0$. However, the probability 
distribution of $\Gnoise$ may not be completely known, and might to some extent
depend on $\ftrue$.

As a regularized inverse of $T$ we consider a class of (non-linear) operators
$S_\lambda: Y \to X$ depending on a regularization parameter $\lambda$ and defined
by
\begin{equation}\label{eq:reginv}
S_\lambda(g) := \argmin_{f \in X} R_\lambda(f) + \frac{1}{2}\|Tf - g\|^2
\end{equation}
where $R_\lambda: X \to \R$ is a regularization functional parameterized by $\lambda$.
The idea is that the reconstruction $\frec = S_\lambda(\gdata)$ might be a good
approximation to $\ftrue$ for suitable choice of the regularization parameter.

Let us assume the existence of a unique minimizer of the optimization 
problem in~\eqref{eq:reginv}. Let us also assume that $R_\lambda$ is convex, and that 
\begin{equation}\label{eq:l1}
R_\lambda(\alpha f) = |\alpha|R_\lambda(f), \quad \forall \alpha \in \R.
\end{equation}

The problem considered here is how to choose the regularization parameter $\lambda$.
Various methods for choosing regularization parameters are of course known.
The method proposed here is motivated by observing the solution of optimization 
problems similar to \eqref{eq:reginv}, but where $f$ is restricted to vary along a line
in $X$. Explicitly, if $f\in X$ and $g \in Y$, define
\begin{equation}\label{eq:alphadef}
\alpha_\lambda(f,g) := \argmin_{\alpha \in \R}  R_\lambda(\alpha f) + \frac{1}{2}\|T(\alpha f) - g\|^2.
\end{equation}
Contrary to \eqref{eq:reginv}, the solution of \eqref{eq:alphadef} can be computed explicitly.

\begin{lemma}\label{lemma1}
Let $f \in X$ and $g \in Y$, and suppose that $R_\lambda$ satisfies~\eqref{eq:l1}. 
If $Tf \ne 0$ or $R_\lambda(f) > 0$ then \eqref{eq:alphadef} has a unique solution given by
\begin{equation}
\alpha_\lambda(f,g) = \left\{
\begin{array}{cc}
\displaystyle \frac{\<Tf, g\> - R_\lambda(f)}{\|Tf\|^2}, & \<Tf, g\> > R_\lambda(f)\\\\
0, & |\<Tf, g\>| \le R_\lambda(f)\\\\
\displaystyle \frac{\<Tf, g\> + R_\lambda(f)}{\|Tf\|^2}, & \<Tf, g\> < -R_\lambda(f).
\end{array}\right.
\end{equation}
If on the other hand $Tf = 0$ and $R_\lambda(f) = 0$, the solution is not unique.
\end{lemma}

The proof is a straightforward computation.

Now we want to look at the random variable $\alpha_\lambda(f,\Gnoise)$. The intuitive idea is
that if $\alpha_\lambda(f, \Gnoise)$ is close to 0 with high probability, this is
an indication that $S_\lambda(\gdata)$ is not heavily influenced by noise. This intuition may or may
not be correct, depending on the nature of the regularization functional and the forward operator.
The precise conditions needed for the idea to be valid are not yet clear.

Let us for any $f \in X$ define 
\begin{equation}\label{eq:sigmadef}
\sigma(f) := \Var[\<Tf, \Gnoise\>]^{1/2}
\end{equation}
and
\begin{equation}\label{eq:sdef}
s_\lambda(f) := \frac{R_\lambda(f)}{\sigma(f)}.
\end{equation}
Lemma~\ref{lemma1} indicates that if $s_\lambda(f) \gg 1$, then $\alpha_\lambda(f,\Gnoise) = 0$ 
with high probability. On the other hand, if $s_\lambda(f) \ll 1$, then $\alpha_\lambda(f,\Gnoise)$
is not strongly affected by the regularization functional. 

The conclusion is that the values of
$s_\lambda(f)$ for different $f \in X$ might serve as quantitative measures of the strength of
the regularization as compared to the noise level. This leads us up to the formulation of
the basic form of the proposed parameter choice rule:
\begin{enumerate}
\item Choose a finite set $F$ of elements in $X$. 
\item For each $f \in F$, choose $\smin(f) \in \R$. The choice of $\smin(f)$ determines how
strongly the inverse problem is regularized. As a rule of thumb, $\smin(f) \ge 5$ corresponds
to very strong regularization, while $\smin(f) \le 1$ corresponds to weak regularization.
\item Choose the smallest regularization parameter $\lambda$ such that $s_\lambda(f) \ge \smin(f)$
for all $f \in F$.
\end{enumerate}

In Section~\ref{sec:TV} a suitable choice of the set $F$ and a model for computing $\smin(f)$ 
is given for the special case of TV regularization. In Section~\ref{sec:ET} a method for estimating
$\sigma(f)$ for ET data is provided. With these ingredients we will then be ready to apply the
parameter choice rule in TV regularized ET.

\section{Application to TV regularization}\label{sec:TV}

Let us specialize the setting as follows. Let $\Omega \subset \R^n$ be a bounded open set
and let $X$ be the space of functions of bounded variation with support in $\overline{\Omega}$, 
with the total variation norm
\begin{equation}\label{eq:tvdef}
\|f\|_{\rm TV} := \sup\left\{\int_{\R^n} f\nabla\cdot h\,dx : h \in C_0^1(\R^n, \R^n),\quad |h(x)| \le 1\right\}
\end{equation}
where $C_0^1(\R^n, \R^n)$ denotes the space of continuously differentiable functions
from $\R^n$ to $\R^n$ with compact support and $\nabla \cdot h$ is the divergence of $h$.
 For continuously differentiable functions $f$ (among others) the total variation is given by
\begin{equation}
\|f\|_{\rm TV} = \int_{\R^n} |\nabla f(x)|\,dx.
\end{equation}
However, the space $X$ contains many non-differentiable functions, including for example the characteristic
functions of certain sets, known as Caccioppoli sets. Caccioppoli sets include all
sets with $C^2$ boundary, see \cite{Gi84} for details.

In TV regularization, the regularization functional is chosen to be $R_\lambda(f) := \lambda\|f\|_{\rm TV}$.
Also, let us assume that $\sigma(f)$ is translation invariant, meaning that if 
$f_1, f_2 \in X$ are related by $f_1(x) = f_2(x-x_0)$ for some $x_0 \in \R^n$, then
$\sigma(f_1) = \sigma(f_2)$. Since $\|f_1\|_{\rm TV} = \|f_2\|_{\rm TV}$, it then follows from the
definition of $s_\lambda$ that $s_\lambda(f_1) = s_\lambda(f_2)$.

\subsection{Choice of the set $F$}

Here we must choose a finite subset $F \subset X$ which is somehow representative of 
all functions in $X$. Lemma~\ref{lemma2} below suggests that it is reasonable to restrict attention
to characteristic functions of Caccioppoli subsets of $\Omega$.

Suppose that $E \subset \Omega$ is a Caccioppoli set, and that $E = E_1 \cup E_2$ with 
$\overline{E}_1 \cap \overline{E}_2 = \emptyset$. Let $f_1$ and $f_2$ be the characteristic functions
of $E_1$ and $E_2$, so $f := f_1 + f_2$ is the characteristic function of $E$. Then 
$\sigma(f) \le \sigma(f_1) + \sigma(f_2)$ and $\|f\|_{\rm TV} = \|f_1\|_{\rm TV} + \|f_2\|_{\rm TV}$
and it follows that 
$$ s_\lambda(f) \ge \lambda \frac{\|f_1\|_{\rm TV} + \|f_2\|_{\rm TV}}{\sigma(f_1) + \sigma(f_2)}
= \frac{s_\lambda(f_1)\sigma(f_1) + s_\lambda(f_2)\sigma(f_2)}{\sigma(f_1) + \sigma(f_2)}
\ge \min \{s_\lambda(f_1), s_\lambda(f_2)\}.$$
Hence, if $f_1$ and $f_2$ are both included in $F$, there is no need to include $f$ unless
$\smin(f) > \min \{\smin(f_1), \smin(f_2)\}$. From this observation it would be tempting to conclude that
only characteristic functions of connected Caccioppoli sets need to be included in $F$. However, 
this is not strictly true, since there are very complicated Caccioppoli sets, which for example
can have uncountably many connected components. Nevertheless, of all functions that are likely to be
treated numerically in practice, only characteristic functions of connected sets need to be 
included in $F$. By similar reasoning, if $\Omega$ has connected complement, only characteristic functions
of sets with connected complement need to be included in $F$.

I further suggest that in many cases it should be reasonable to choose $F$ as a set of characteristic
functions of balls of different sizes. The precise arguments for this and the conditions under which
they are valid remain to be clarified. 

Suppose we make the choice to let $F$ consist of characteristic functions of balls. We would then
choose a finite set $D = \{d_1, \ldots, d_k\}$ of diameters of these balls. By the assumption
of translation invariance, it is sufficient to include one ball of each diameter in $F$,
so we have $F = \{f_d: d \in D\}$ where $f_d$ denotes the characteristic function of an 
arbitrary ball of diameter $d$.

To conclude this section, we state and prove the lemma which was used above to motivate the 
restriction to characteristic functions of Caccioppoli sets.

\begin{lemma}\label{lemma2}
Suppose there exists a constant $C$ such that $\sigma(f) \le C\|f\|_{\rm TV}$ and
$\sigma(f) \le C\|f\|_{L^\infty}$ for all $f \in X$.
If $s_\lambda(f) \ge s_0 > 0$ for every $f\in X$ which is the characteristic function of 
a Caccioppoli subset of $\Omega$, then the same inequality holds for every $f \in X$.
\end{lemma}

\begin{remark}
From the hypothesis of the lemma it is trivially true that $s_\lambda(f) \ge \lambda/C$.
However, the constant $C$ could a priori be very large. The point of the lemma is that
if a better estimate holds for characteristic functions of sets, then the same 
estimate necessarily holds for all functions of bounded variation.
\end{remark}

\begin{proof}
The idea of the proof is simply that an arbitrary function is a superposition of
characteristic functions defined by its level sets.

Note first that it is sufficient to prove the statement for positive functions.
For assuming that this has been done, an arbitrary function can be written as
$f = f_+ - f_-$ where $f_+$ and $f_-$ are positive and $\|f\|_{\rm TV}
= \|f_+\|_{\rm TV} + \|f_-\|_{\rm TV}$. Since $\sigma(f) \le \sigma(f_+) + \sigma(f_-)$
it follows that
$$s_\lambda(f) \ge \lambda\frac{\|f_+\|_{\rm TV} + \|f_-\|_{\rm TV}}{\sigma(f_+) + \sigma(f_-)}
= \frac{s_\lambda(f_+)\sigma(f_+) + s_\lambda(f_-)\sigma(f_-)}{\sigma(f_+) + \sigma(f_-)}\ge s_0.$$

So suppose that $f$ is a positive function. For $t > 0$ define $\chi_t$ to be
the characteristic function of the set $\{x: f(x) \ge t\}$. By the coarea formula
\cite[Theorem 1.23]{Gi84}, it holds that
\begin{equation}
\|f\|_{\rm TV} = \int_0^\infty \|\chi_t\|_{\rm TV}\,dt.
\end{equation}

If it can be shown that
\begin{equation}\label{eq:sigmaineq}
\sigma(f) \le \int_0^\infty \sigma(\chi_t)\,dt
\end{equation}
the conclusion follows, since then
$$s_\lambda(f) \ge \frac{\lambda \int_0^\infty \|\chi_t\|_{\rm TV}\,dt}{\int_0^\infty \sigma(\chi_t)\,dt}
= \frac{\int_0^\infty s_\lambda(\chi_t)\sigma(\chi_t)\,dt}{\int_0^\infty \sigma(\chi_t)\,dt}
\ge s_0.$$

Let $0 = t_0 < t_1 < t_2 < \cdots < t_m$ be a sequence of positive numbers, and define
$\delta_i = t_i - t_{i-1}$. For each $i = 1, \ldots, m$ there is some $\tau_i \in [t_{i-1}, t_i]$
such that
$$ \delta_i \sigma(\chi_{\tau_i}) \le \int_{t_{i-1}}^{t_i} \sigma(\chi_t)\,dt. $$
If
$$ f_1 = \sum_{i=1}^m \delta_i \chi_{\tau_i} $$
then it holds that
\begin{equation}\label{eq:sigmaineq2}
\sigma(f_1) \le \sum_{i=1}^m \delta_i \sigma(\chi_{\tau_i})
\le \int_0^{t_m} \sigma(\chi_t)\,dt.
\end{equation}
If we define $f_2(x) = \max\{0, f(x) - t_m\}$, it follows from the coarea formula 
and the assumptions of the lemma that
$$ \sigma(f_2) \le C\|f_2\|_{\rm TV} = C\int_{t_m}^\infty \|\chi_t\|_{\rm TV}\, dt. $$
Finally, we have that
$$ \sigma(f - f_1 - f_2) \le C\|f - f_1 - f_2\|_{L^\infty} \le C\max_{1\le i\le m} \delta_i.$$
This shows that $\sigma(f_2)$ and $\sigma(f - f_1 - f_2))$ can be made arbitrarily small by making
$m$ and $t_m$ large. Since $\sigma(f) \le \sigma(f_1) + \sigma(f_2) + \sigma(f - f_1 - f_2))$,
\eqref{eq:sigmaineq} follows from \eqref{eq:sigmaineq2} and this completes the proof.
\end{proof}
 
\subsection{Model for computing $\smin(f)$}

Here I will provide a model for computing $\smin(f)$. Let $f_d \in F$ be a function whose 
support is a ball of diameter $d$, for example a characteristic function as suggested in the
previous section. Let us assume that we are given 
a real number $a \ge 0$ and want to avoid that $|\frec(x)| > a$ in regions where 
$\ftrue(x) = 0$. Let us also assume that the probability distribution of $\<Tf, \Gnoise\>$
is to a good approximation Gaussian.

Consider balls of diameter $d$ in $\Omega$. The maximum number $N_d$ of such disjoint balls in $\Omega$
is approximately $N_d \approx |\Omega|d^{-n}$. Let $f_d^1, \ldots , f_d^{N_d}$ be translations of
$f_d$ with support in these disjoint balls. A heuristic argument suggests that we should look at 
the probability that $|\alpha_\lambda(f_d^j, \Gnoise)|\cdot \|f_d^j\|_{L^\infty} > a$. By 
Lemma~\ref{lemma1}, that probability is (assuming that $\<Tf_d^j, \Gnoise\>$ is Gaussian)
\begin{equation}
\erfc\left(\frac{1}{\sqrt{2}}\left(s_\lambda(f_d) + \frac{a\|Tf_d\|^2}
{\sigma(f_d)\|f_d\|_{L^\infty}}\right)\right).
\end{equation}
Let us, rather arbitrarily, choose the regularization parameter so that the expected number 
of $j \in \{1, \ldots, N_d\}$ with $|\alpha_\lambda(f_d^j, \Gnoise)| > a$ is not more than 1.
This is equivalent to the inequality
\begin{equation}\label{eq:smin}
\begin{split}
s_\lambda(f_d) \ge \smin(f_d) &= \sqrt{2}\erfcinv\left(\frac{1}{N_d}\right) - \frac{a\|Tf_d\|^2}{\sigma(f_d)\|f_d\|_{L^\infty}}\\
&\approx \sqrt{2}\erfcinv\left(\frac{d^n}{|\Omega|}\right) - \frac{a\|Tf_d\|^2}{\sigma(f_d)\|f_d\|_{L^\infty}}.
\end{split}
\end{equation}
The formula is easily modified if some other restriction on the expected number of $j$ with 
$|\alpha_\lambda(f_d^j, \Gnoise)| > a$ is desired.

\section{A numerical example from electron tomography}\label{sec:ET}

In this section I will apply the proposed parameter choice method to a numerical example
from electron tomography (ET). First I give a brief description of this inverse problem.
For details, see for example \cite{FaOk08}.

Electron tomography is a method using a transmission electron microscope to construct three-dimensional
models of biological macromolecules and similar structures. The data collected consists of a series 
of images, a {\em tilt series}, with the specimen tilted in different angles. Hence the data space can be
decomposed as a direct sum $Y = Y_1 \oplus \cdots \oplus Y_m$, where $Y_j$ corresponds to the
$j$th image in the tilt series, and each vector $g \in Y$ can be written as $g = g_1 + \cdots + g_m$
where $g_j \in Y_j$.
The forward operator $T$ has a corresponding decomposition into components $T_j: X \to Y_j$.

A commonly used approximation of the forward operator is that each $T_j$ consists of a parallel
beam transform in a direction depending on $j$, followed by convolution with a point spread
function. The noise in the data comes mainly from the stochastic nature of the detection of the electrons,
and has a Poisson distribution. Due to the necessity of using a very low electron dose 
in each image, the noise level is usually very high.

In this model, the assumption of translation invariance made in Section~\ref{sec:TV}
is valid. Since $\<Tf, \Gnoise\>$ is composed of noise from all the
images, which can reasonably be assumed to be independent, the assumption that its probability
distribution is Gaussian seems plausible according to the central limit theorem.

\subsection{Estimation of $\sigma(f)$}

In order to apply the proposed parameter choice method, it is necessary to estimate 
$\sigma(f)$ for a given $f \in X$. For tomographic data of the type encountered in ET, this
estimate can be made directly from the data set, given the following very reasonable assumptions.

\begin{enumerate}
\item The noise components in separate images are uncorrelated. \label{assum1}
\item The noise components in different parts of the same image are at most weakly 
correlated.\label{assum2}
\item The probability distribution of $\<T_j f, \Gnoise_j\>$ is invariant under translation 
of $f$.\label{assum3}
\item In each image, the signal to noise ratio is much lower than 1.\label{assum4}
\end{enumerate}

Now, by assumption~\ref{assum1} we have that
\begin{equation}
\Var[\<Tf, \Gnoise\>] = \sum_{j=1}^m \Var[\<T_jf, \Gnoise_j\>].
\end{equation}
To estimate the right hand side of this equality, take a number of random translations
$f_1, \ldots, f_l$ of $f$ (this requires that the support of $f$ is small compared to $\Omega$).
I claim that $\Var[\<T_jf, \Gnoise_j\>]$ can be approximated by a sample variance
\begin{equation}\label{eq:samplevar}
\Var[\<T_jf, \Gnoise_j\>] \approx \frac{1}{l-1}\sum_{i=1}^l \left(\<T_jf_i, \gnoise_j\>
- \frac{1}{l}\sum_{i'=1}^l \<T_jf_{i'}, \gnoise_j\>\right)^2.
\end{equation}
This is justified by assumptions~\ref{assum2} and \ref{assum3}. Finally, assumption~\ref{assum4}
justifies that we can replace $\gnoise$ by $\gdata$ in \eqref{eq:samplevar}. Combining these
steps leads to the following approximation:
\begin{equation}
\sigma(f)^2 \approx \sum_{j=1}^m \frac{1}{l-1}\sum_{i=1}^l \left(\<T_jf_i, \gdata_j\>
- \frac{1}{l}\sum_{i'=1}^l \<T_jf_{i'}, \gdata_j\>\right)^2.
\end{equation}

\subsection{Numerical results}

The numerical results presented in this section were obtained by approximately solving
the minimization problem \eqref{eq:reginv} with $T$ the forward operator from electron
tomography described above
and $g$ a real or simulated data set. The regularization functional is an 
approximation of the total variation norm, defined as follows. Let 
$x_{i,j,k}, (i,j,k)\in I \subset \Z^3$ be a 
rectangular lattice of points at which the function $f$ is sampled. For all 
$(i,j,k)\notin I$ we take $f(x_{i,j,k})$ to be 0. Let $\bar{I}$ be the subset of $\Z^3$
consisting of $I$ together with all points adjacent to $I$. Let 
\begin{gather*}
D_1^+f(x_{i,j,k}) := f(x_{i+1,j,k}) - f(x_{i,j,k})\\
D_2^+f(x_{i,j,k}) := f(x_{i,j+1,k}) - f(x_{i,j,k})\\
D_3^+f(x_{i,j,k}) := f(x_{i,j,k+1}) - f(x_{i,j,k})
\end{gather*}
be discrete partial derivatives of $f$ in the forward direction, and similarly let
\begin{gather*}
D_1^-f(x_{i,j,k}) := f(x_{i,j,k}) - f(x_{i-1,j,k})\\
D_2^-f(x_{i,j,k}) := f(x_{i,j,k}) - f(x_{i,j-1,k})\\
D_3^-f(x_{i,j,k}) := f(x_{i,j,k}) - f(x_{i,j,k-1})
\end{gather*}
be discrete partial derivatives in the backward direction for all $(i,j,k)\in \bar{I}$. 
Now we define
\begin{equation}
R_\lambda(f) = \lambda \sum_{\bar{I}}\left(\beta^2 
+ \frac{1}{2}\sum_{l=1}^3 \Bigl((D_l^+ f(x_{i,j,k}))^2 + (D_l^- f(x_{i,j,k}))^2\Bigr)\right)^{1/2}.
\end{equation}
The parameter $\beta$ is included in order to make the regularization functional smooth,
which was necessary for the minimization algorithm used. (For alternative optimization algorithms
which do not require this approximation, see for example \cite{Hi04}.) It was set to a small positive value,
$\beta = 3\cdot10^{-4}$. This is well below the level where changes in $\beta$ do not seem to have
any noticeable effect on the solution. However, in the application of the parameter choice 
rule, $\beta$ was set to 0, so that the condition \eqref{eq:l1} is exactly satisfied.

The approximate solution of the minimization problem \eqref{eq:reginv} was computed by 
iteratively searching for the minimum in 2-dimensional subspaces of $X$, where each 
subspace is spanned by the gradient of the objective functional and a vector in the direction of
the previous update. This method seems to be considerably faster than a gradient descent
method minimizing over a 1-dimensional subspace in each iteration. The iteration was 
continued until no appreciable change occured even after many iterations. The number of 
iterations used was in most cases between 500 and 1500, with a fairly good approximation 
of the end result occuring within 100 iterations.

\subsubsection{Simulated data set}

\begin{figure}[htp]
\begin{center}
\includegraphics[scale=0.3]{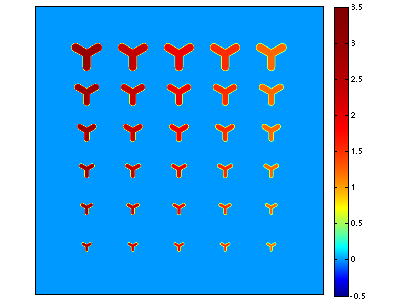}
\includegraphics[scale=0.3]{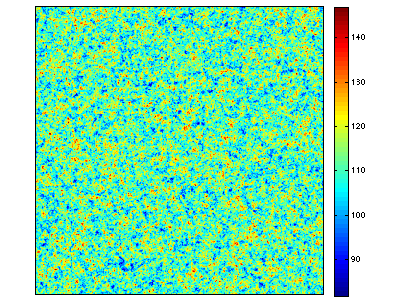}
\end{center}
\caption{{\em Left}\/: A section through the phantom.
{\em Right}\/: The central projection in the data set.\label{fig:1}}
\end{figure}

The first numerical example is reconstructed from a simulated data set. The phantom used in the
simulation contains 30 Y-shaped objects of varying size and contrast. From this phantom an ET
data set was simulated, consisting of 121 projections. The specimen was tilted about a single
axis, with the tilt angle ranging from $-60^\circ$ to $60^\circ$. The simulated electron dose
was 15.7 electrons per pixel on average over the tilt series. A section through the phantom and
the central projection in the data set are shown in figure~\ref{fig:1}.

\begin{figure}[htp]
\begin{center}
\includegraphics[scale=0.3]{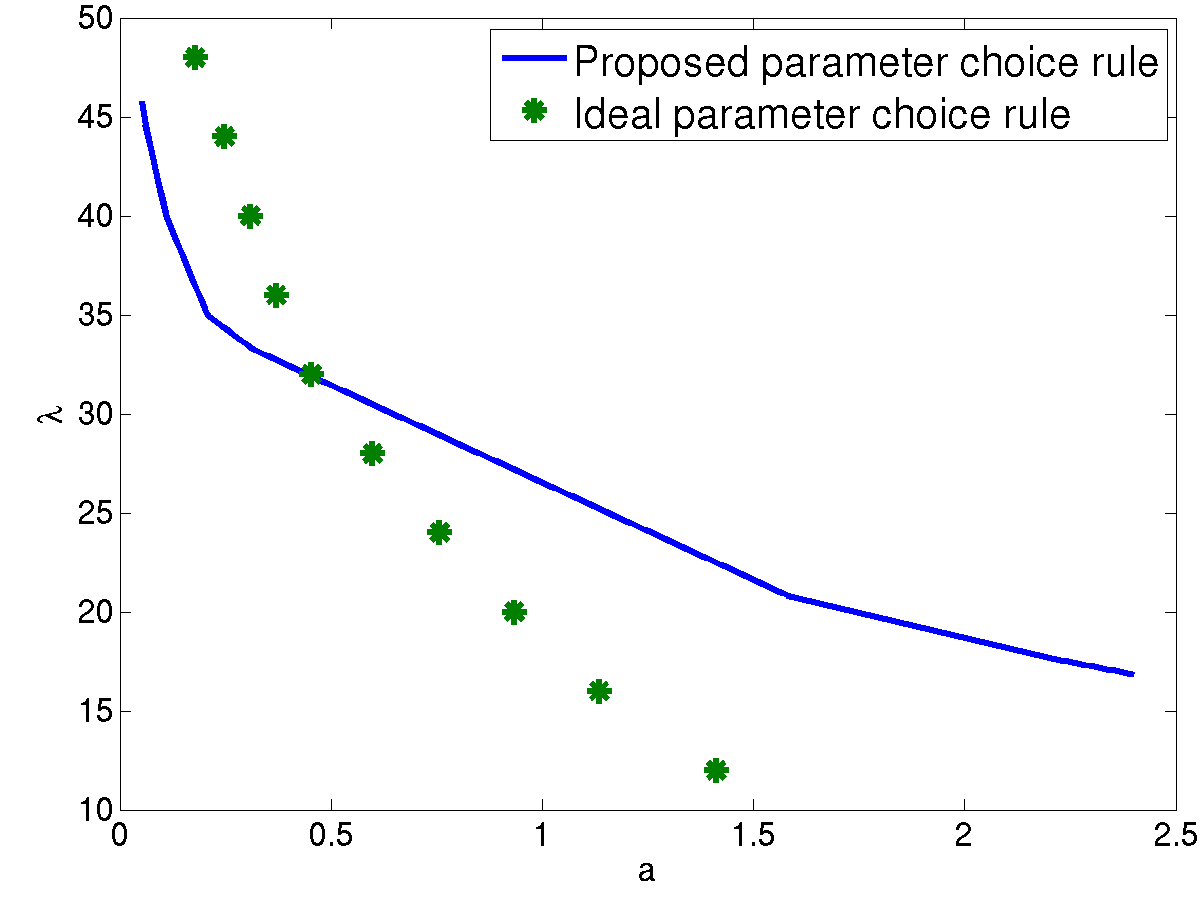}
\includegraphics[scale=0.3]{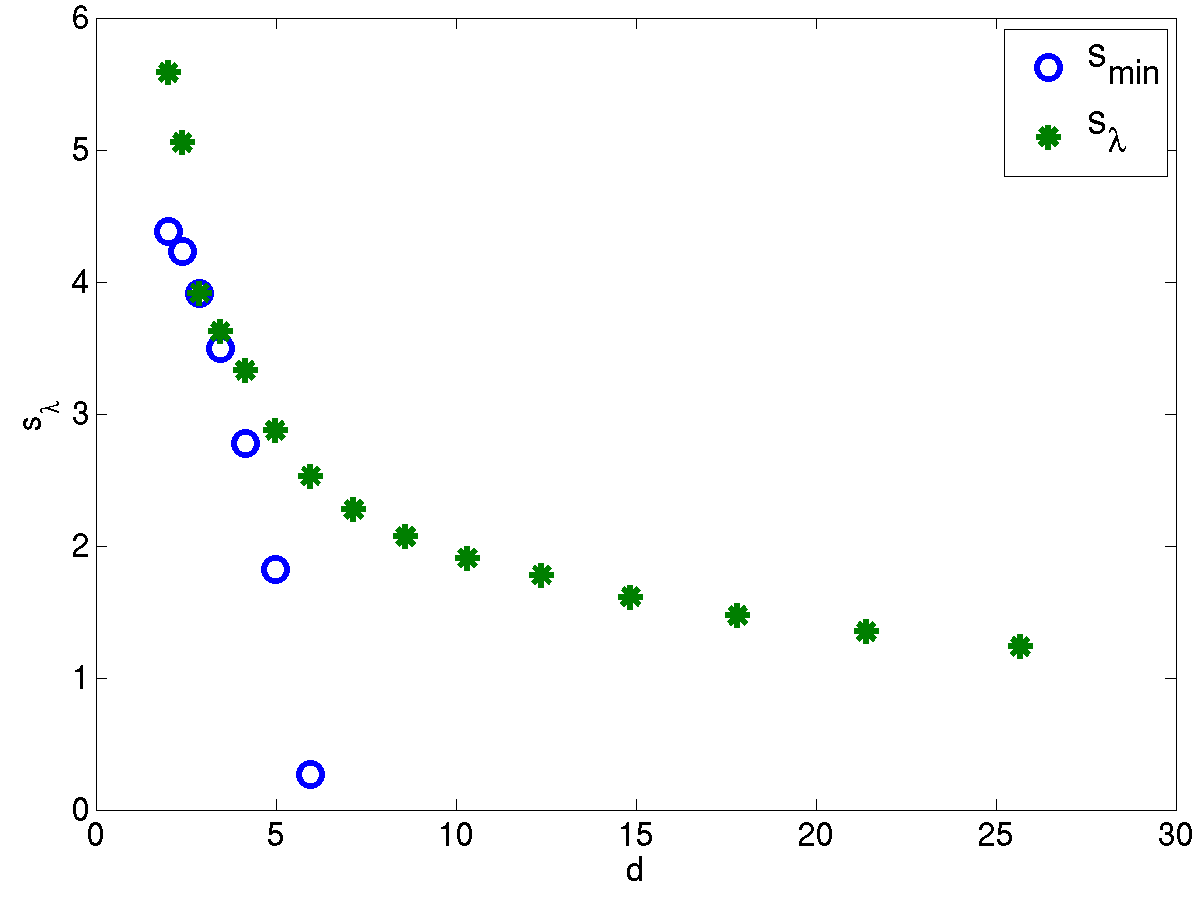}
\end{center}
\caption{{\em Left}\/: The solid line shows the dependence of the regularization parameter $\lambda$ on the 
threshold parameter $a$ when it is chosen according to the proposed parameter choice rule. The graph 
is composed of line segments since a discrete set of diameters $d$ were used in the computation. 
This graph is compared to the ideal parameter choice rule for this problem (stars), which 
can be determined by comparing a set of reconstructions with different regularization parameters 
to the true solution (see below).
{\em Right}\/: The dependence of $\smin(f_d)$ (circles) and $s_\lambda(f_d)$
(stars) on $d$ for $a = 0.5$ and the corresponding $\lambda = 31.4$. For $d$ greater than 
approximately 6, $\smin(f_d)$ drops below 0, and does not impose any restriciton on $\lambda$.
\label{fig:2}}
\end{figure}

Next the parameter choice rule from the previous sections was applied to the data set. Note that
the only input needed is the data set, a model for the forward operator and the threshold parameter $a$. 
A reasonable choice for the threshold parameter, given the overall levels of contrast in the
phantom, seems to be $a = 0.5$. With this choice we should hope that the reconstructed objects
are well above the noise level, even if the contrast is somewhat reduced by the regularization.
Figure~\ref{fig:2} shows the dependence of $\lambda$ on $a$, and the dependence of $\smin(f_d)$ 
and $s_\lambda(f_d)$ on $f_d$ for the choice of $\lambda$ corresponding to $a = 0.5$. The
diameter $d$ is measured in voxel units.

\begin{figure}[htp]
\begin{center}
\includegraphics[scale=0.3]{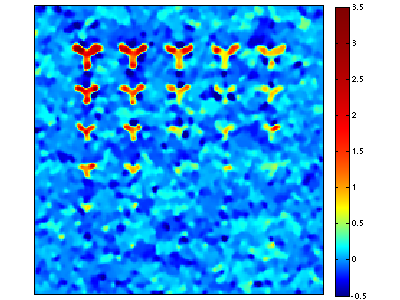}
\includegraphics[scale=0.3]{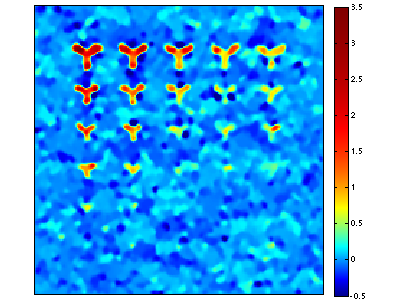}
\end{center}
\caption{{\em Left}\/: Section through reconstruction with $\lambda = 12$.
{\em Right}\/: Section through reconstruction with $\lambda = 16$.\label{fig:3}}
\end{figure}

\begin{figure}[htp]
\begin{center}
\includegraphics[scale=0.3]{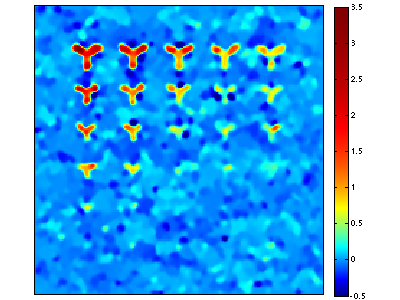}
\includegraphics[scale=0.3]{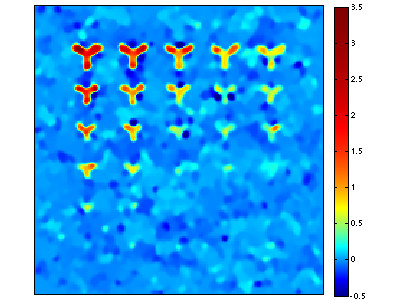}
\end{center}
\caption{{\em Left}\/: Section through reconstruction with $\lambda = 20$.
{\em Right}\/: Section through reconstruction with $\lambda = 24$.\label{fig:4}}
\end{figure}

\begin{figure}[htp]
\begin{center}
\includegraphics[scale=0.3]{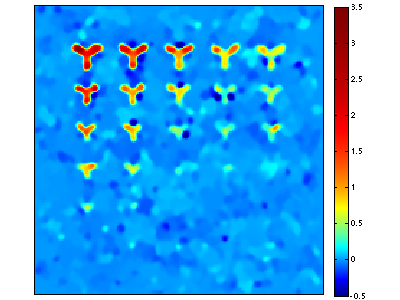}
\includegraphics[scale=0.3]{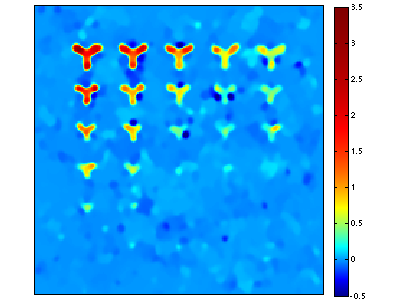}
\end{center}
\caption{{\em Left}\/: Section through reconstruction with $\lambda = 28$.
{\em Right}\/: Section through reconstruction with $\lambda = 32$.
These are in the range obtained by applying the proposed parameter selection rule, depending
on the threshold parameter.\label{fig:5}}
\end{figure}

\begin{figure}[htp]
\begin{center}
\includegraphics[scale=0.3]{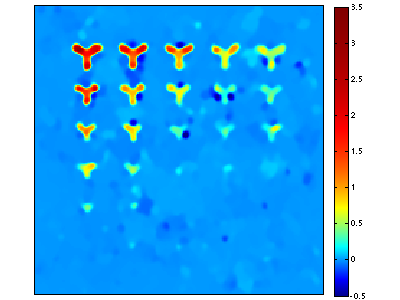}
\includegraphics[scale=0.3]{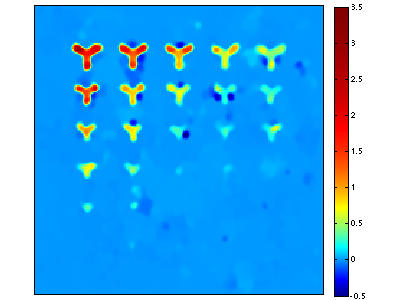}
\end{center}
\caption{{\em Left}\/: Section through reconstruction with $\lambda = 36$.
{\em Right}\/: Section through reconstruction with $\lambda = 40$.\label{fig:6}}
\end{figure}

\begin{figure}[htp]
\begin{center}
\includegraphics[scale=0.3]{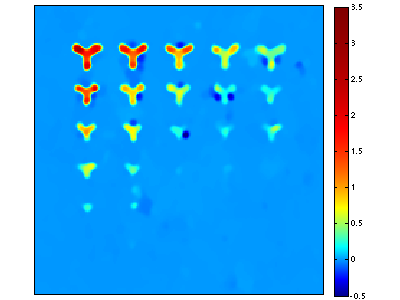}
\includegraphics[scale=0.3]{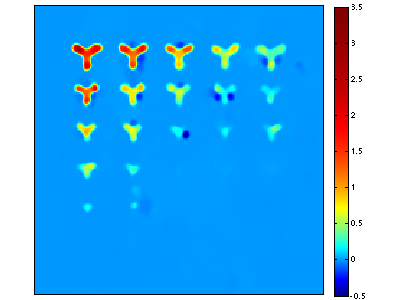}
\end{center}
\caption{{\em Left}\/: Section through reconstruction with $\lambda = 44$.
{\em Right}\/: Section through reconstruction with $\lambda = 48$.\label{fig:7}}
\end{figure}

If the proposed parameter choice rule is applicable, $\lambda \approx 30$ should be a suitable 
choice of the regularization parameter. To test if this is the case, a series of TV regularized
reconstructions were computed with the regularization parameter ranging from 12 to 48. The size of
the reconstructions is $200 \times 200 \times 100$ voxels. A section
through each of the reconstructions is shown in figures~\ref{fig:3}--\ref{fig:7}.
Which one of these reconstructions would be considered optimal is of course strongly dependent
on the type of further analysis it is intended for.

To further investigate the amount of undesirable noise in the reconstructions, the following
analysis was applied. For a given threshold $a \ge 0$ and a given reconstruction
$\frec$, the set $\{x \in \Omega: \frec(x) > a\}$ was computed and decomposed into
its connected components. (In the discrete setting, a voxel was considered to be connected 
to each of its 8 nearest neighbors.) A connected component was classified as a 
true hit if it has nonempty intersection with some object in the phantom, otherwise
it was classified as a false hit. Hence, the number of true and false hits can be counted,
where the count depends both on the regularization parameter $\lambda$ and the threshold
$a$. Only one true hit was counted for each of the objects in the phantom, so the
number of true hits can never exceed 30. Table~\ref{table:1} shows how the number of true
and false hits depend on the regularization parameter when $a = 0.5$. For $\lambda$ below
a certain level, the number of false hits increases dramatically, and there is an 
obvious risk of misinterpretation. In this case, the proposed parameter choice
rule makes a surprisingly accurate prediction of the point where false hits start to 
occur.

\begin{table}[htp]
\begin{center}
\begin{tabular}{|l|cccccccccc|}
\hline
$\lambda$ & 12 & 16 & 20 & 24 & 28 & 32 & 36 & 40 & 44 & 48\\
\hline
True hits & 25 & 24 & 22 & 19 & 19 & 16 & 14 & 13 & 13 & 12\\
\hline
False hits & 471 & 98 & 28 & 5 & 1 & 0 & 1 & 0 & 0 & 0\\
\hline
\end{tabular}
\end{center}
\caption{The number of true and false hits at threshold level $0.5$ as a function
of the regularization parameter $\lambda$.\label{table:1}}
\end{table}

Another way to look at false hits is as follows. In a given reconstruction we can 
compute the smallest threshold $a$ which does not give rise to any false hits. This defines
a relation between $a$ and $\lambda$ which can be considered as an ideal parameter
choice rule. 
The only obvious way to determine this ideal parameter choice rule is
to compute a set of reconstructions with different regularization parameters
and compare them to the true solution $\ftrue$, in order to classify reconstructed objects as
true or false hits. This, of course, impossible in real life problems, primarily because $\ftrue$ 
is not known. However, with simulated data, the ideal parameter choice rule can be
compared to a practically applicable parameter choice rule. 

Such a comparison is shown in figure~\ref{fig:1} above. The plot shows that, while there is 
certainly a discrepancy, the proposed parameter choice rule provides at least a rough 
approximation of the ideal parameter choice rule for this particular problem. The regularization
parameter selected by the proposed rule tends to be too small when $a$ is small and too 
large when $a$ is large. The reason for this trend is not yet clear.

\subsubsection{Real data set of TMV specimen}

The second example uses a real ET data set of a specimen containing Tobacco Mosaic Virus (TMV).
The TMV is a long and fairly rigid cylindrical object, approximately \unit{18}{\nano\meter}
in diameter. The tilt series contains 61 projections, with the tilt angle varying in the range
$-60^\circ$ to $60^\circ$. The electron dose was 64.5 electrons per pixel on average over the
tilt series.

\begin{figure}[htp]
\begin{center}
\includegraphics[scale=0.3]{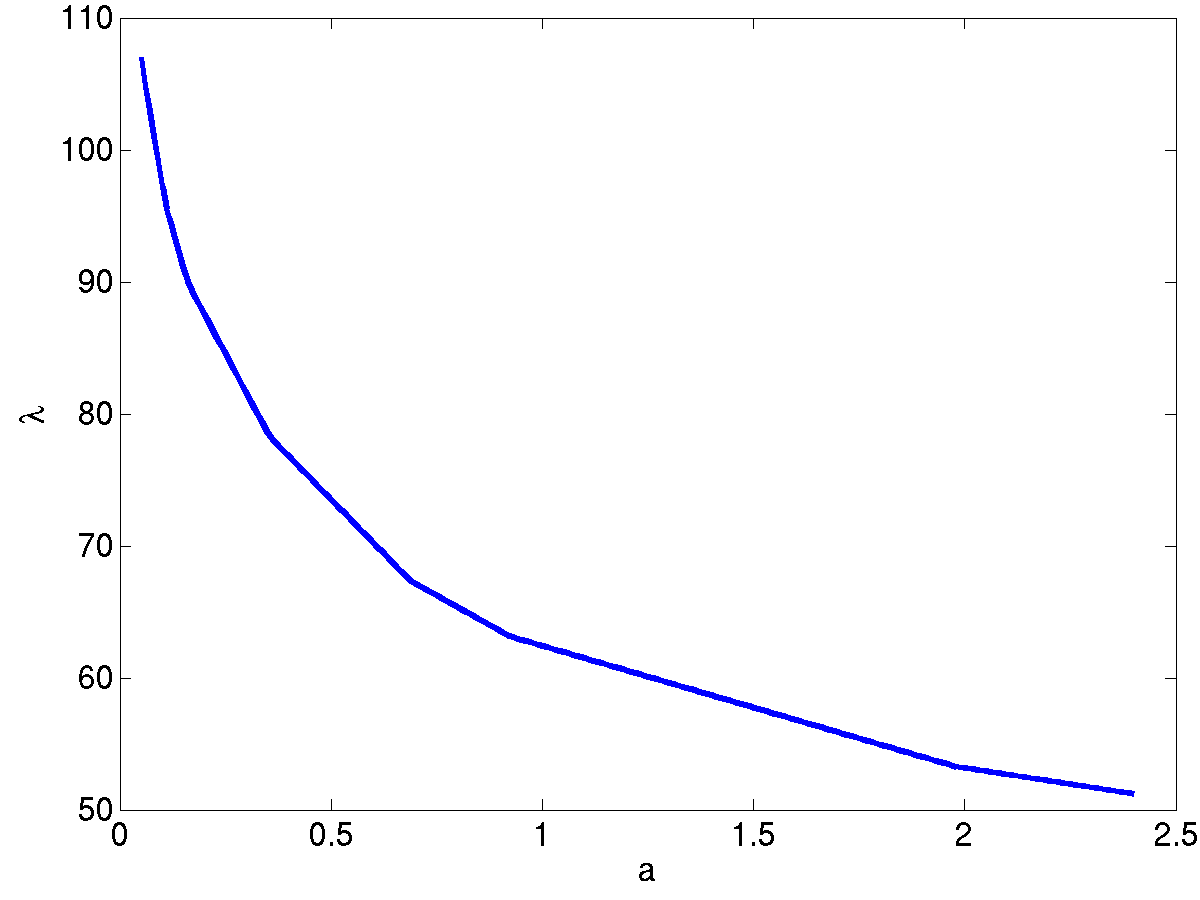}
\includegraphics[scale=0.3]{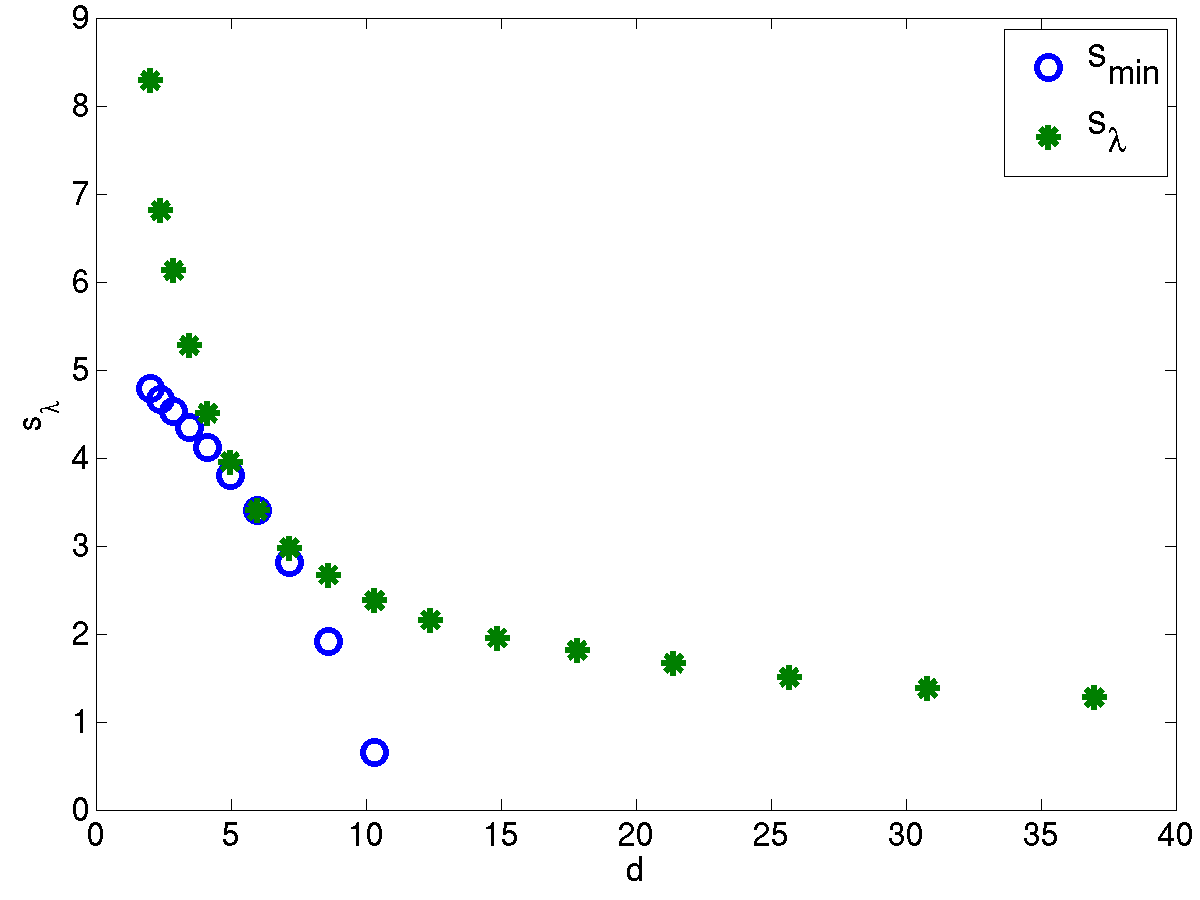}
\end{center}
\caption{{\em Left}\/: The dependence of the regularization parameter $\lambda$ on the threshold
parameter $a$ when it is chosen according to the proposed parameter choice rule.
{\em Right}\/: The dependence of $\smin(f_d)$ (circles) and $s_\lambda(f_d)$
(stars) on $d$ for $a = 0.5$ and the corresponding $\lambda = 73.5$.\label{fig:8}}
\end{figure}

Application of the parameter choice rule to the data set yields the dependence of the 
regularization parameter $\lambda$ on the threshold $a$ shown in figure~\ref{fig:8}. This 
suggests that $\lambda \approx 70$ should be a suitable choice.

\begin{figure}[htp]
\begin{center}
\includegraphics[scale=0.3]{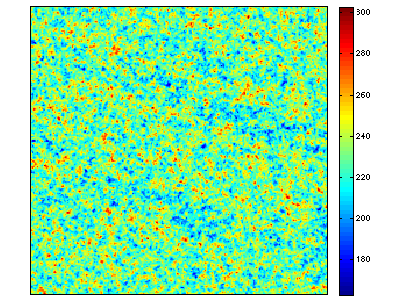}
\includegraphics[scale=0.3]{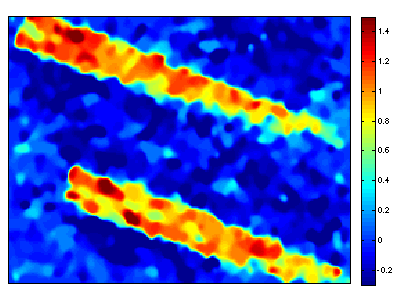}
\end{center}
\caption{{\em Left}\/: The central projection in the data set.
{\em Right}\/: Section through reconstruction with $\lambda = 40$.\label{fig:9}}
\end{figure}

\begin{figure}[htp]
\begin{center}
\includegraphics[scale=0.3]{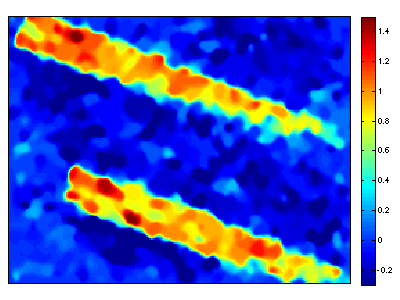}
\includegraphics[scale=0.3]{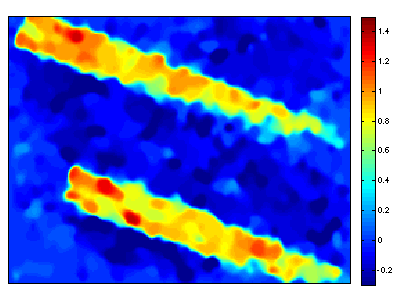}
\end{center}
\caption{{\em Left}\/: Section through reconstruction with $\lambda = 50$.
{\em Right}\/: Section through reconstruction with $\lambda = 60$.\label{fig:10}}
\end{figure}

\begin{figure}[htp]
\begin{center}
\includegraphics[scale=0.3]{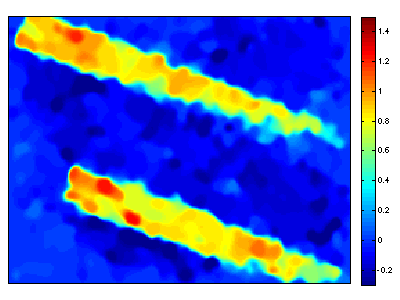}
\includegraphics[scale=0.3]{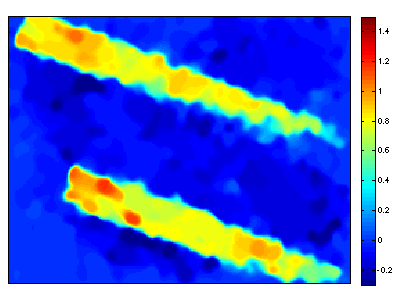}
\end{center}
\caption{{\em Left}\/: Section through reconstruction with $\lambda = 70$.
{\em Right}\/: Section through reconstruction with $\lambda = 80$.\label{fig:11}}
\end{figure}

Sections through reconstructions with a range of regularization parameters are shown in 
figures~\ref{fig:9}--\ref{fig:11}. The reconstructed volume contains two TMV's, which are
clearly visible in the images. 

When working with real data it is of course not possible to know for certain which 
features in a reconstruction correspond to real objects in the specimen. In this case, 
all that we can be reasonably sure about is that the large elongated objects in the 
reconstructions indicate the presence of TMV's in the specimen. Smaller objects
might be due to noise, but could also represent some contamination in the specimen.
Table~\ref{table:2} lists the number of small objects above the threshold level 
$0.5$ in the reconstructed volume for different values of the regularization parameter.

\begin{table}[htp]
\begin{center}
\begin{tabular}{|l|ccccc|}
\hline
Regularization parameter $\lambda$ & 40 & 50 & 60 & 70 & 80 \\
\hline
Number of small objects & 129 & 40 & 17 & 5 & 1 \\
\hline
\end{tabular}
\end{center}
\caption{The number of small objects, not classified as TMV's, at threshold level
$a = 0.5$ as a function of the regularization parameter $\lambda$.\label{table:2}}
\end{table}

The number of small objects in the reconstruction grows rapidly when the value of the
regularization parameter goes below the value provided by the parameter choice rule.
Although it is impossible to know for sure if these small objects are due only to noise,
the conclusion must be that this might very well be the case. Hence, no significance
should be attached to these objects when the reconstruction is interpreted.

\section{Discussion}\label{sec:discussion}

The numerical examples presented suggest that the parameter choice rule outlined in this
paper might be useful for inverse problems of the type encountered in ET.
To further investigate under what circumstances the method would be useful, it is 
indispensable to gain a better understanding of the underlying mathematics. Once this 
has been done it should be possible to devise a pertinent set of numerical test cases.

It might certainly be argued that in these examples, some interpretations of details 
could be more easily made from the reconstructions with a smaller regularization parameter
than the one indicated by the proposed parameter choice rule. This suggests that perhaps 
\eqref{eq:smin} should be modified so that the computed value of $\smin$ is somewhat 
smaller. This could be done by relaxing the restriction on the expected 
number of $j$ with $|\alpha_\lambda(f_d^j, \Gnoise)| > a$, which was in the derivation 
above, rather arbitrarily, chosen as 1. However, if $\smin$ is decreased, both the 
heuristic argument and the numerical examples indicate that one should expect to 
have false objects appearing in the reconstructions, and the interpretations must be
made with this in mind.

In the application to ET data exemplified in this paper, the proposed parameter choice rule
is heuristic, or error free, the magnitude of the noise being estimated directly from the
data. It is well known that such heuristic parameter choice rules, 
can not provide convergent regularization methods in the strict sense, that is, methods that converge
to the true solution when the noise level approaches zero. Such heuristic parameter choice rules
might nevertheless be very useful in practice, when the noise level certainly does not
approach zero. As mentioned in the introduction, one of the motivations for this
new parameter choice rule is the difficulties associated with very high noise levels.
A variant of the proposed parameter choice rule would be to estimate the quantities 
$\sigma(f)$ from some a priori estimate of the noise level. Used in this way, the parameter choice 
could provide a convergent regularization method.

In the derivations leading up to the proposed parameter choice rule, particularly in the
estimation of $\sigma(f)$, rather restrictive, albeit very realistic, assumptions were made
on the statistical properties of the noise. It is worth noting that the necessity of making
such assumptions, explicitly or implicitly, is inherent to the inverse problem in ET; 
without them no reconstruction would be possible at all. A malicious demon, if allowed to 
choose the noise with only a restriction on its norm, could easily hide every trace of the 
signal we are trying to recover. 

One advantage of the suggested parameter choice rule, as compared for example to the L-curve 
method, is its computational efficiency.
In the numerical examples considered here, the application of the parameter choice rule
has a lower computational cost than one single iteration of the algorithm subsequently
used to compute the regularized solution.

The considerations leading up to the parameter choice rule can also be employed in an
alternative way. Suppose a regularization parameter has been chosen by some other 
method and a reconstruction has been computed. In the interpretation of the reconstruction
it is then desirable to know which features are reliable, and which are likely to be
an effect of noise. A computation of the quantity $s_\lambda(f)$, where $f$ is a certain
feature in the reconstruction, can then be used to estimate how likely it is for such 
features to occur solely as a consequence of random noise.

\paragraph{Acknowledgements} The author wants to thank Ozan \"Oktem for encouragement and support, and
Sidec AB for supplying ET data.

\bibliographystyle{plain}
\bibliography{invprob}

\begin{thebibliography}{1}

\bibitem{ChEsPaYi06}
T.~Chan, S.~Esedoglu, F.~Park, and A.~Yip.
\newblock Total variation image restoration: overview and recent developments.
\newblock In {\em Handbook of mathematical models in computer vision}, pages
  17--31. Springer, New York, 2006.

\bibitem{EnHaNe96}
Heinz~W. Engl, Martin Hanke, and Andreas Neubauer.
\newblock {\em Regularization of inverse problems}, volume 375 of {\em
  Mathematics and its Applications}.
\newblock Kluwer Academic Publishers Group, Dordrecht, 1996.

\bibitem{FaOk08}
D.~Fanelli and O.~{\"O}ktem.
\newblock Electron tomography: a short overview with an emphasis on the
  absorption potential model for the forward problem.
\newblock {\em Inverse Problems}, 24(1):013001, 2008.

\bibitem{Gi84}
Enrico Giusti.
\newblock {\em Minimal surfaces and functions of bounded variation}, volume~80
  of {\em Monographs in Mathematics}.
\newblock Birkh\"auser Verlag, Basel, 1984.

\bibitem{Hi04}
M.~Hinterm{\"u}ller and K.~Kunisch.
\newblock Total bounded variation regularization as a bilaterally constrained
  optimization problem.
\newblock {\em SIAM J. Appl. Math.}, 64(4):1311--1333 (electronic), 2004.

\end{thebibliography}

\end{document}